\documentclass[preprint,12pt]{elsarticle}

\usepackage{amsmath,amssymb,amsfonts,enumerate,amsthm}
\usepackage{dsfont,mathrsfs,euscript}
\usepackage{blkarray,float,hhline}

\usepackage{pgf,tikz,pgfplots}
\pgfplotsset{compat=1.14}
\usetikzlibrary{arrows}
\usetikzlibrary{matrix,decorations.pathreplacing}
\usepackage[colorlinks=true,citecolor=blue,linkcolor=red,urlcolor=blue]{hyperref}

\newtheorem{theorem}{Theorem}[section]

\newtheorem{conjecture}[theorem]{Conjecture}

\theoremstyle{definition}

\newtheorem{example}[theorem]{Example}

\theoremstyle{remark}

\journal{}

\begin{document}

\begin{frontmatter}



\title{A counterexample to a conjecture of Wang and Hou on signed graphs}


\author{Asghar Bahmani}

\address{Department of Mathematics and Computer Science, Amirkabir University of Technology, P.O. Box: 15875-4413, Tehran, Iran}
\ead{asghar.bahmani@gmail.com,asghar.bahmani@aut.ac.ir}
 
\begin{abstract}
We give a counterexample to a conjecture of Wang and Hou  related with the sum of the $k$ largest Laplacian eigenvalues of signed graphs.
\end{abstract}


\end{frontmatter}


\section{Introduction}

Let $G$ be a simple graph. A signed graph $G_{\sigma}$ is obtained from $G$ by a sign function $\sigma: E(G)\rightarrow \{+1,-1\}$. The adjacency matrix $A(G_{\sigma})$ is a $(0,\pm 1)$-matrix where
\[
A(G_{\sigma})_{ij}=\begin{cases}
\sigma(ij) & ij\in E(G),\\ 
0 & \text{othewise.}\\
\end{cases}
\]
For a positive integer $n$, we denote by $[n]$, the set $\{1,\ldots,n\}$. The Laplacian matrix of a simple graph $G$, denoted by $L_{G}$, is $D(G)-A(G)$, where $D(G)$ is a diagonal matrix of the degree sequence of $G$ and $A(G)$ is the adjacency matrix of $G$.
Similar to the $L_{G}$, the Laplacian matrix of a signed graph $G_{\sigma}$ is defined as $D(G)-A(G_{\sigma})$ and is denoted by $L_{G_{\sigma}}$. The following well-known conjecture states an upper bounds for the sum of the largest eigenvalues of $L_{G}$.

\begin{conjecture}[Brouwer's Conjecture]\cite{bh}
	Let $k,n\in \mathbb{N}$ and $k\in[n]$. For any simple graph $G$ of order $n$, if $\mu_1(L_G)\geq\cdots\geq\mu_n(L_G)$ are the Laplacian eigenvalues of $G$, then we have $\sum_{i=1}^{k}\mu_i(L_G)\leq e(G)+\binom{k+1}{2}$. 
\end{conjecture}

In the recent paper of Wang and Hou, "Wang D., Hou Y., On the sum of Laplacian eigenvalues of a signed graph, Linear Algebra and its Applications, 555(2018), 39-52", they conjectured an upper bound for the sum of the $k$ largest Laplacian eigenvalues of signed graphs.

\begin{conjecture}\label{whc}\cite{WH}
	Let $k,n\in \mathbb{N}$ and $k\in[n]$. For any signed graph $G_{\sigma}$, if $\mu_1(L_{G_{\sigma}})\geq\cdots\geq\mu_n(L_{G_{\sigma}})$ are the Laplacian eigenvalues of $G_{\sigma}$, then we have $\sum_{i=1}^{k}\mu_i(L_{G_{\sigma}})\leq e(G)+\binom{k+1}{2}+1$. 
\end{conjecture}

We considered Conjecture \ref{whc} for graphs of order $7$ and $8$ and by using computer-assisted computations to signing the complete graphs of order $7$ and $8$, we found some counterexamples for $k=4$ and $k=5$.

\begin{example}\label{74e}
For the graph $G_{\sigma}$ as shown below, Figure \ref{74}, we have
	\[L_{G_{\sigma}}=
\begin{pmatrix}
6&	1&	1&	-1&	1&	-1&	1\\
1&	6&	1&	1&	-1&	1&	1\\
1&	1&	6&	1&	1&	-1&	-1\\
-1&	1&	1&	6&	1&	1&	-1\\
1&	-1&	1&	1&	6&	1&	1\\
-1&	1&	-1&	1&	1&	6&	1\\
1&	1&	-1&	-1&	1&	1&	6
\end{pmatrix} \,\, \text{and}\]
\[
 (\mu_1(L_{G_{\sigma}}),\cdots,\mu_7(L_{G_{\sigma}}))\simeq (8.7015,8.2360,8.2360,7.,3.7639,3.7639,2.2984).\] So, $\sum_{i=1}^{4}\mu_i(L_{G_{\sigma}})\simeq 32.1735\nleq e(G)+\binom{5}{2}+1=32$. 
\end{example}

\begin{figure}[H]\label{74}
	\centering
	\includegraphics[scale=0.4]{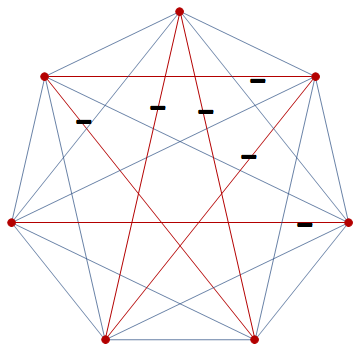}
	\caption{The graph $G_{\sigma}$ in Example \ref{74e}.}
\end{figure}

\begin{example}\label{85e}
	For the graph $H_{\sigma}$ as shown below, Figure \ref{85}, we have 
	\[L_{H_{\sigma}}=
	\begin{pmatrix}
	7&	1&	1&	1&	-1&	1&	1&	1\\
	1&	7&	1&	1&	1&	-1&	-1&	1\\
	1&	1&	7&	1&	1&	-1&	1&	-1\\
	1&	1&	1&	7&	1&	1&	-1&	-1\\
	-1&	1&	1&	1&	7&	1&	1&	1\\
	1&	-1&	-1&	1&	1&	7&	1&	1\\
	1&	-1&	1&	-1&	1&	1&	7&	1\\
	1&	1&	-1	&-1&	1&	1&	1&	7
	\end{pmatrix}\,\, \text{and} \]
\[
	(\mu_1(L_{H_{\sigma}}),\cdots,\mu_8(L_{H_{\sigma}}))\simeq(10.6056,10.,8.,8.,8.,4.,4.,3.3944).\] 
	So, $\sum_{i=1}^{5}\mu_i(L_{H_{\sigma}})\simeq 44.6056\nleq e(H)+\binom{6}{2}+1=44$. 
\end{example}

\begin{figure}[H]\label{85}
	\centering
	\includegraphics[scale=.4]{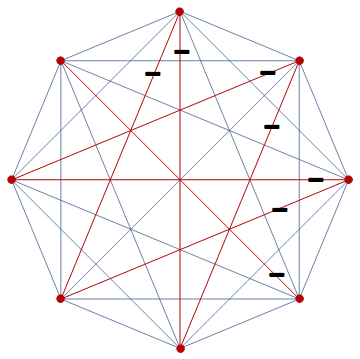}
	\caption{The graph $H_{\sigma}$ in Example \ref{85e}.}
\end{figure}

\section*{Acknowledgments}
The author is indebted to Iran National Science Foundation (INSF) for supporting this research under grant number 95005902.
 




\end{document}